# SUFFICIENT BURN-IN FOR GIBBS SAMPLERS FOR A HIERARCHICAL RANDOM EFFECTS MODEL

By Galin L. Jones and James P. Hobert[1]

*University of Minnesota and University of Florida*

We consider Gibbs and block Gibbs samplers for a Bayesian hierarchical version of the one-way random effects model. Drift and minorization conditions are established for the underlying Markov chains. The drift and minorization are used in conjunction with results from J. S. Rosenthal [*J. Amer. Statist. Assoc.* **90** (1995) 558–566] and G. O. Roberts and R. L. Tweedie [*Stochastic Process. Appl.* **80** (1999) 211–229] to construct analytical upper bounds on the distance to stationarity. These lead to upper bounds on the amount of burn-in that is required to get the chain within a prespecified (total variation) distance of the stationary distribution. The results are illustrated with a numerical example.

**1. Introduction.** We consider a Bayesian hierarchical version of the standard normal theory one-way random effects model. The posterior density for this model is intractable in the sense that the integrals required for making inferences cannot be computed in closed form. Hobert and Geyer (1998) analyzed a Gibbs sampler and a block Gibbs sampler for this problem and showed that the Markov chains underlying these algorithms converge to the stationary (i.e., posterior) distribution at a geometric rate. However, Hobert and Geyer stopped short of constructing analytical upper bounds on the total variation distance to stationarity. In this article, we construct such upper bounds and this leads to a method for determining a sufficient *burn-in*.

Our results are useful from a practical standpoint because they obviate troublesome, ad hoc convergence diagnostics [Cowles and Carlin (1996) and

Received January 2002; revised December 2002.
[1]Supported in part by NSF Grant DMS-00-72827.

Cowles, Roberts and Rosenthal (1999)]. More important, however, we believe that this is the first analysis of a *practically relevant* Gibbs sampler on a continuous state space that provides viable burn-ins. By practically relevant, we mean that the stationary distribution is complex enough that independent and identically distributed (i.i.d.) sampling is not straightforward. We note that the Gibbs samplers analyzed by Hobert (2001) and Rosenthal (1995a, 1996) are not practically relevant since i.i.d. samples can be drawn from the corresponding stationary distributions using simple, sequential sampling schemes [Jones (2001) and Marchev and Hobert (2004)]. Some notation is now introduced that will allow for a more detailed overview.

Let $X = \{X_i, i = 0, 1, \ldots\}$ be a discrete time, time homogeneous Markov chain that is irreducible, aperiodic and positive Harris recurrent. Let $P^n(x, \cdot)$ be the probability measure corresponding to the random variable $X_n$ conditional on starting the chain at $X_0 = x$; that is, $P^n$ is the $n$-step Markov transition kernel. Let $\pi(\cdot)$ be the invariant probability measure of the chain and let $\|\cdot\|$ denote the total variation norm. Formally, the issue of burn-in can be described as follows. Given a starting value $x_0$ and an arbitrary $\varepsilon > 0$, can we find an $n^* = n^*(x_0, \varepsilon)$ such that

$$\|P^{n^*}(x_0, \cdot) - \pi(\cdot)\| < \varepsilon? \tag{1}$$

If the answer is "yes," then, since the left-hand side of (1) is nonincreasing in the number of iterations, the distribution of $X_k$ is within $\varepsilon$ of $\pi$ for all $k \geq n^*$. Because we are not demanding that $n^*$ be the smallest value for which (1) holds, it is possible that the chain actually gets within $\varepsilon$ of stationarity in much fewer than $n^*$ iterations. For this reason, we call $n^*$ a *sufficient burn-in*.

Several authors [see, e.g., Meyn and Tweedie (1994), Rosenthal (1995a), Cowles and Rosenthal (1998), Roberts and Tweedie (1999) and Douc, Moulines and Rosenthal (2002)] have recently provided results that allow one to calculate $n^*$ when $X$ is *geometrically ergodic*. However, to use these results one must establish both a *drift condition* and an associated *minorization condition* for $X$. [For an accessible treatment of these concepts, see Jones and Hobert (2001).] In this article we establish drift and minorization for the Gibbs samplers analyzed by Hobert and Geyer (1998). These conditions are used in conjunction with the theorems of Rosenthal (1995a) and Roberts and Tweedie (1999) to construct formulas that can be used to calculate $n^*$.

The rest of the article is organized as follows. The model and algorithms are described in Section 2. Section 3 contains important background material on general state space Markov chain theory as well as statements of the theorems of Rosenthal (1995a) and Roberts and Tweedie (1999). This section also contains a new *conversion lemma* that provides a connection between the two different types of drift used in these theorems. We establish drift and minorization for the block Gibbs sampler in Section 4 and the same is done for the Gibbs sampler in Section 5. In Section 6 the results are illustrated and Rosenthal's theorem is compared with the theorem of Roberts and Tweedie. Section 7 contains some concluding remarks.

**2. The model and the Gibbs samplers.** Consider the following Bayesian version of the standard normal theory one-way random effects model. First, conditional on $\theta = (\theta_1, \ldots, \theta_K)^T$ and $\lambda_e$ the data $Y_{ij}$ are independent with

$$Y_{ij}|\theta, \lambda_e \sim \mathrm{N}(\theta_i, \lambda_e^{-1}),$$

where $i = 1, \ldots, K$ and $j = 1, \ldots, m_i$. At the second stage, conditional on $\mu$ and $\lambda_\theta$, $\theta_1, \ldots, \theta_K$ and $\lambda_e$ are independent with

$$\theta_i|\mu, \lambda_\theta \sim \mathrm{N}(\mu, \lambda_\theta^{-1}) \quad \text{and} \quad \lambda_e \sim \mathrm{Gamma}(a_2, b_2),$$

where $a_2$ and $b_2$ are known positive constants. [We say $W \sim \mathrm{Gamma}(\alpha, \beta)$ if its density is proportional to $w^{\alpha-1} e^{-w\beta} I(w > 0)$.] Finally, at the third stage $\mu$ and $\lambda_\theta$ are assumed independent with

$$\mu \sim \mathrm{N}(m_0, s_0^{-1}) \quad \text{and} \quad \lambda_\theta \sim \mathrm{Gamma}(a_1, b_1),$$

where $m_0, s_0, a_1$ and $b_1$ are known constants; all but $m_0$ are assumed to be positive so that all of the priors are proper. The posterior density of this hierarchical model is characterized by

$$(2) \qquad \pi_h(\theta, \mu, \lambda|y) \propto f(y|\theta, \lambda_e) f(\theta|\mu, \lambda_\theta) f(\lambda_e) f(\mu) f(\lambda_\theta),$$

where $\lambda = (\lambda_\theta, \lambda_e)^T$, $y$ is a vector containing all of the data, and $f$ denotes a generic density. [We will often abuse notation and use $\pi_h$ to denote the probability distribution associated with the density in (2).] Expectations with respect to $\pi_h$ are typically ratios of intractable integrals, the numerators of which can have dimension as high as $K+3$ [Jones and Hobert (2001)]. Thus, to make inferences using $\pi_h$, we must resort to (possibly) high dimensional numerical integration, analytical approximations or Monte Carlo and Markov chain Monte Carlo techniques.

In their seminal article on the Gibbs sampler, Gelfand and Smith (1990) used the balanced version of this model (in which $m_i \equiv m$) as an example. [See also Gelfand, Hills, Racine-Poon and Smith (1990) and Rosenthal (1995b).] Each iteration of the standard, fixed-scan Gibbs sampler consists of updating all of the $K + 3$ variables in the same predetermined order. The *full conditionals* required for this Gibbs sampler are now reported. Define

$$v_1(\theta, \mu) = \sum_{i=1}^{K} (\theta_i - \mu)^2,$$

$$v_2(\theta) = \sum_{i=1}^{K} m_i (\theta_i - \bar{y}_i)^2 \quad \text{and} \quad \mathrm{SSE} = \sum_{i,j} (y_{ij} - \bar{y}_i)^2,$$

where $\bar{y}_i = m_i^{-1} \sum_{j=1}^{m_i} y_{ij}$. The full conditionals for the variance components are

$$\lambda_\theta | \theta, \mu, \lambda_e, y \sim \text{Gamma}\left(\frac{K}{2} + a_1, \frac{v_1(\theta, \mu)}{2} + b_1\right) \tag{3}$$

and

$$\lambda_e | \theta, \mu, \lambda_\theta, y \sim \text{Gamma}\left(\frac{M}{2} + a_2, \frac{v_2(\theta) + \text{SSE}}{2} + b_2\right), \tag{4}$$

where $M = \sum_i m_i$. Letting $\theta_{-i} = (\theta_1, \ldots, \theta_{i-1}, \theta_{i+1}, \ldots, \theta_K)^T$ and $\bar{\theta} = K^{-1} \times \sum_i \theta_i$, the remaining full conditionals are

$$\theta_i | \theta_{-i}, \mu, \lambda_\theta, \lambda_e, y \sim \text{N}\left(\frac{\lambda_\theta \mu + m_i \lambda_e \bar{y}_i}{\lambda_\theta + m_i \lambda_e}, \frac{1}{\lambda_\theta + m_i \lambda_e}\right)$$

for $i = 1, \ldots, K$ and

$$\mu | \theta, \lambda_\theta, \lambda_e, y \sim \text{N}\left(\frac{s_0 m_0 + K \lambda_\theta \bar{\theta}}{s_0 + K \lambda_\theta}, \frac{1}{s_0 + K \lambda_\theta}\right).$$

We consider the fixed-scan Gibbs sampler that updates $\mu$, then the $\theta_i$'s, then $\lambda_\theta$ and $\lambda_e$. Since the $\theta_i$'s are conditionally independent given $(\mu, \lambda)$, the order in which they are updated is irrelevant. The same is true of $\lambda_\theta$ and $\lambda_e$ since these two random variables are conditionally independent given $(\theta, \mu)$. If we write a one-step transition as $(\mu', \theta', \lambda') \to (\mu, \theta, \lambda)$, then the Markov transition density (MTD) of our Gibbs sampler is given by

$$k(\mu, \theta, \lambda | \mu', \theta', \lambda') = f(\mu | \theta', \lambda'_\theta, \lambda'_e, y) \left[\prod_{i=1}^{K} f(\theta_i | \theta_{-i}, \mu, \lambda'_\theta, \lambda'_e, y)\right]$$

$$\times f(\lambda_\theta | \theta, \mu, \lambda'_e, y) f(\lambda_e | \theta, \mu, \lambda_\theta, y).$$

Hobert and Geyer (1998) considered this same update order. We note here that, in general, Gibbs samplers with different update orders correspond to different Markov chains. However, two chains whose update orders are cyclic permutations of one another converge at the same rate.

As an alternative to the standard Gibbs sampler, Hobert and Geyer (1998) introduced the more efficient *block* Gibbs sampler in which all of the components of $\xi = (\theta_1, \ldots, \theta_K, \mu)^T$ are updated simultaneously. These authors showed that $\xi | \lambda, y \sim \text{N}(\xi^*, V)$ and gave formulas for $\xi^* = \xi^*(\lambda, y)$ and $V = V(\lambda, y)$. Because we will make extensive use of these formulas, they are restated in Appendix A. One iteration of the block Gibbs sampler consists of updating $\lambda_\theta$, $\lambda_e$ and $\xi$ in some order. Due to the conditional independence of $\lambda_\theta$ and $\lambda_e$, the block Gibbs sampler is effectively a two-variable Gibbs sampler or *data augmentation* algorithm [Tanner and Wong (1987)], the two components being $\xi$ and $\lambda$. We choose to update $\lambda$ first because, as we will see later, updating the most complicated distribution last typically simplifies the calculations required to establish drift and minorization conditions. If we write a one-step transition as $(\lambda', \xi') \to (\lambda, \xi)$, then the corresponding MTD is given by

$$
\begin{aligned}
(5) \quad k(\lambda, \xi | \lambda', \xi') &= f(\lambda | \xi', y) \, f(\xi | \lambda, y) \\
&= f(\lambda_\theta | \xi', y) \, f(\lambda_e | \xi', y) \, f(\xi | \lambda_\theta, \lambda_e, y).
\end{aligned}
$$

Hobert and Geyer (1998) considered the opposite update order because they were not attempting to simultaneously establish drift *and* minorization. Note, however, that our update order is just a cyclic permutation of the order used by Hobert and Geyer.

A proper formulation of the burn-in problem requires some concepts and notation from Markov chain theory. These are provided in the following section. More general accounts of this material can be found in Nummelin (1984), Meyn and Tweedie (1993) and Tierney (1994).

**3. Markov chain background.** Let $\mathcal{X} \subset \mathbb{R}^p$ for $p \geq 1$ and let $\mathcal{B}$ denote the associated Borel $\sigma$-algebra. Suppose that $X = \{X_i, i = 0, 1, \dots\}$ is a discrete time, time homogeneous Markov chain with state space $\mathcal{X}$ and Markov transition kernel $P$; that is, for $x \in \mathcal{X}$ and $A \in \mathcal{B}$, $P(x, A) = \Pr(X_{i+1} \in A | X_i = x)$. Also, for $n = 1, 2, 3, \dots$, let $P^n$ denote the $n$-step transition kernel, that is, $P^n(x, A) = \Pr(X_{i+n} \in A | X_i = x)$ so, in particular, $P \equiv P^1$. Note that $P^n(x, \cdot)$ is the probability measure of the random variable $X_n$ conditional on starting the chain at $X_0 = x$.

Let $\nu$ be a measure on $\mathcal{B}$. We will say that the Markov chain $X$ satisfies assumption $(\mathcal{A})$ if it is $\nu$-irreducible, aperiodic and positive Harris recurrent with invariant probability measure $\pi(\cdot)$. It is straightforward to show that the Gibbs samplers described in the previous section satisfy assumption $(\mathcal{A})$ with $\nu$ equal to Lebesgue measure. Under assumption $(\mathcal{A})$, for every $x \in \mathcal{X}$ we have

$$\|P^n(x, \cdot) - \pi(\cdot)\| \downarrow 0 \qquad \text{as } n \to \infty,$$

where $\|P^n(x, \cdot) - \pi(\cdot)\| := \sup_{A \in \mathcal{B}} |P^n(x, A) - \pi(A)|$ is the total variation distance between $P^n$ and $\pi$. The chain $X$ is called *geometrically ergodic* if it satisfies assumption $(\mathcal{A})$ and, in addition, there exist a constant $0 < t < 1$ and a function $g : \mathcal{X} \mapsto [0, \infty)$ such that, for any $x \in \mathcal{X}$,

$$(6) \qquad \|P^n(x, \cdot) - \pi(\cdot)\| \leq g(x) t^n$$

for $n = 1, 2, \dots$. It has recently been demonstrated that establishing drift and minorization conditions for $X$ verifies geometric ergodicity (the existence of $g$ and $t$) and yields an upper bound on the right-hand side of (6). See Jones and Hobert (2001) for an expository look at this theory. In this paper, we will focus on the results due to Rosenthal (1995a) and Roberts and Tweedie (1999). Slightly simplified versions of these results follow.

THEOREM 3.1 [Rosenthal (1995a)]. *Let $X$ be a Markov chain satisfying assumption ($\mathcal{A}$). Suppose $X$ satisfies the following drift condition. For some function $V : \mathcal{X} \mapsto [0, \infty)$, some $0 < \gamma < 1$ and some $b < \infty$,*

$$(7) \qquad E[V(X_{i+1})|X_i = x] \le \gamma V(x) + b \qquad \forall\, x \in \mathcal{X}.$$

*Let $C = \{x \in \mathcal{X} : V(x) \le d_\mathrm{R}\}$, where $d_\mathrm{R} > 2b/(1-\gamma)$ and suppose that $X$ satisfies the following minorization condition. For some probability measure $Q$ on $\mathcal{B}$ and some $\varepsilon > 0$,*

$$(8) \qquad P(x, \cdot) \ge \varepsilon Q(\cdot) \qquad \forall\, x \in C.$$

*Let $X_0 = x_0$ and define two constants as follows:*

$$\alpha = \frac{1 + d_\mathrm{R}}{1 + 2b + \gamma d_\mathrm{R}} \quad \text{and} \quad U = 1 + 2(\gamma d_\mathrm{R} + b).$$

*Then, for any $0 < r < 1$,*

$$\|P^n(x_0, \cdot) - \pi(\cdot)\| \le (1-\varepsilon)^{rn} + \left(\frac{U^r}{\alpha^{1-r}}\right)^n \left(1 + \frac{b}{1-\gamma} + V(x_0)\right).$$

THEOREM 3.2 [Roberts and Tweedie (1999, 2001)]. *Let $X$ be a Markov chain satisfying assumption ($\mathcal{A}$). Suppose $X$ satisfies the following drift condition. For some function $W : \mathcal{X} \mapsto [1, \infty)$, some $0 < \rho < 1$ and some $L < \infty$,*

$$(9) \qquad E[W(X_{i+1})|X_i = x] \le \rho W(x) + L I_S(x) \qquad \forall\, x \in \mathcal{X},$$

*where $S = \{x \in \mathcal{X} : W(x) \le d_\mathrm{RT}\}$ and*

$$d_\mathrm{RT} \ge \frac{L}{1-\rho} - 1.$$

*Suppose further that $X$ satisfies the following minorization condition. For some probability measure $Q$ on $\mathcal{B}$ and some $\varepsilon > 0$,*

$$(10) \qquad P(x, \cdot) \ge \varepsilon Q(\cdot) \qquad \forall\, x \in S.$$

*Let $X_0 = x_0$ and define some constants as follows:*

$$\kappa = \rho + \frac{L}{1 + d_\mathrm{RT}}, \qquad J = \frac{(\kappa d_\mathrm{RT} - \varepsilon)(1 + d_\mathrm{RT}) + L d_\mathrm{RT}}{(1 + d_\mathrm{RT})\kappa},$$

$$\zeta = \frac{\log[(1/2)(L/(1-\rho) + w(x_0))]}{\log(\kappa^{-1})}, \qquad \eta = \frac{\log[(1-\varepsilon)^{-1} J]}{\log(\kappa^{-1})},$$

$$\beta_\mathrm{RT} = \exp\left[\frac{\log \kappa \log(1-\varepsilon)}{\log J - \log(1-\varepsilon)}\right].$$

*Then if $J \ge 1$ and $n' = k - \zeta > \eta(1-\varepsilon)/\varepsilon$, we have, for any $1 \le \beta < \beta_\mathrm{RT}$,*

$$(11) \quad \|P^k(x_0, \cdot) - \pi(\cdot)\| < \left[1 - \frac{\beta(1-\varepsilon)}{(1 + \eta/n')^{1/\eta}}\right]\left(1 + \frac{n'}{\eta}\right)\left(1 + \frac{\eta}{n'}\right)^{n'/\eta} \beta^{-n'}.$$

REMARK 3.1. The version of Theorem 3.2 in Roberts and Tweedie (1999) relies on their Theorem 5.2, whose proof contains an error. Using Roberts and Tweedie's (1999) notation, suppose $V : \mathcal{X} \mapsto [1, \infty)$, $d > 0$, $C = \{x \in \mathcal{X} : V(x) \leq d\}$ and $h(x,y) = (V(x) + V(y))/2$. Roberts and Tweedie (1999) claim that

$$h(x,y) \geq (1+d) I_{[C \times C]^c}(x,y),$$

which is false and, in fact, all that we can claim is that

$$h(x,y) \geq \frac{1+d}{2} I_{[C \times C]^c}(x,y).$$

We have accounted for this error in our statement of Theorem 3.2 and we are grateful to an anonymous referee for bringing the error to our attention.

REMARK 3.2. Roberts and Tweedie (1999) provide a different bound for the case $J < 1$ but, since we do not use it in our application (see Section 6), it is not stated here.

REMARK 3.3. Roberts and Tweedie (1999) show that the right-hand side of (11) is approximately minimized when $\beta = \beta_{\text{RT}}/(1 + \eta/n')^{1/\eta}$.

REMARK 3.4. It is well known [see, e.g., Meyn and Tweedie (1993), Chapter 15] that (7) and (8) together [or (9) and (10) together] imply that $X$ is geometrically ergodic. See Jones and Hobert (2001) for an heuristic explanation.

In our experience it is often easier to establish a Rosenthal-type drift condition than a Roberts-and-Tweedie-type drift condition. The following new result provides a useful connection between these two versions of drift.

LEMMA 3.1. *Let $X$ be a Markov chain satisfying assumption $(\mathcal{A})$. Suppose there exist $V : \mathcal{X} \mapsto [0, \infty)$, $\gamma \in (0,1)$ and $b < \infty$ such that*

(12) $$E[V(X_{n+1})|X_n = x] \leq \gamma V(x) + b \qquad \forall x \in \mathcal{X}.$$

*Set $W(x) = 1 + V(x)$. Then, for any $a > 0$,*

(13) $$E[W(X_{n+1})|X_n = x] \leq \rho W(x) + L I_C(x) \qquad \forall x \in \mathcal{X},$$

*where $\rho = (a + \gamma)/(a + 1)$, $L = b + (1 - \gamma)$ and*

$$C = \left\{ x \in \mathcal{X} : W(x) \leq \frac{(a+1)L}{a(1-\rho)} \right\}.$$

PROOF. Clearly, (12) implies that

$$E[W(X_{i+1})|X_i = x] \leq \gamma W(x) + b + (1-\gamma) = \gamma W(x) + L \qquad \forall\, x \in \mathcal{X}.$$

Set $\Delta W(x) = E[W(X_{n+1})|X_n = x] - W(x)$ and $\beta = (1-\gamma)/(a+1)$. Then

$$E[W(X_{n+1})|X_n = x] \leq [1 - (a+1)\beta]W(x) + L$$

or, equivalently,

$$\Delta W(x) \leq -\beta W(x) - a\beta W(x) + L$$

for all $x \in \mathcal{X}$. If $x \notin C$, then

$$W(x) > \frac{(a+1)L}{a(1-\rho)} > \frac{(a+1)L}{a(1-\gamma)} = \frac{L}{a\beta}.$$

Now write $W(x) = \frac{L}{a\beta} + s(x)$, where $s(x) > 0$. Then

$$\Delta W(x) \leq -\beta W(x) - a\beta \left[\frac{L}{a\beta} + s(x)\right] + L$$
$$= -\beta W(x) - a\beta s(x)$$
$$\leq -\beta W(x).$$

If, on the other hand, $x \in C$, then

$$\Delta W(x) \leq -\beta W(x) - a\beta W(x) + L$$
$$\leq -\beta W(x) + L.$$

Now putting these together gives

$$E[W(X_{n+1})|X_n = x] \leq (1-\beta)W(x) + L I_C$$
$$= \rho W(x) + L I_C. \qquad \square$$

REMARK 3.5. Since

$$\frac{(a+1)L}{a(1-\rho)} \geq \frac{L}{1-\rho} - 1,$$

(13) constitutes a drift condition of the form (9). Therefore, if we can establish (12) as well as a minorization condition on the set $C$, it will be as straightforward to apply Theorem 3.2 as it is to apply Theorem 3.1. Indeed, this is the approach we take with our Gibbs samplers. Moreover, we use $a = 1$ in our application since $\frac{(a+1)L}{a(1-\rho)}$ is minimized at this value.

While the Gibbs sampler is easier to implement than the block Gibbs sampler, it is actually harder to analyze because it is effectively a three-variable Gibbs sampler as opposed to the block Gibbs sampler, which is effectively a two-variable Gibbs sampler. Thus, we begin with block Gibbs.

**4. Drift and minorization for the block Gibbs sampler.** Drift conditions of the form (7) are established for the unbalanced and balanced cases in Sections 4.1 and 4.2, respectively. A minorization condition that works for both cases is established in Section 4.3. Throughout this section we assume that $m' = \min\{m_1, m_2, \ldots, m_K\} \geq 2$ and that $K \geq 3$.

4.1. *Drift*: *unbalanced case.* Define two constants as follows:

$$\delta_1 = \frac{1}{2a_1 + K - 2} \quad \text{and} \quad \delta_2 = \frac{1}{2a_2 + M - 2}.$$

Also define $\delta_3 = (K+1)\delta_2$ and $\delta_4 = \delta_2 \sum_{i=1}^{K} m_i^{-1}$. Our assumptions about $K$ and $m'$ guarantee that $0 < \delta_i < 1$ for $i = 1, 2, 3, 4$. Set $\delta = \max\{\delta_1, \delta_3\}$. Also, let $\Delta$ denote the length of the convex hull of the set $\{\bar{y}_1, \bar{y}_2, \ldots, \bar{y}_K, m_0\}$ and define

$$c_1 = \frac{2b_1}{2a_1 + K - 2} \quad \text{and} \quad c_2 = \frac{2b_2 + \text{SSE}}{2a_2 + M - 2}.$$

PROPOSITION 4.1. *Fix $\gamma \in (\delta, 1)$ and let $\phi_1$ and $\phi_2$ be positive numbers such that $\frac{\phi_1 \delta_4}{\phi_2} + \delta < \gamma$. Define the drift function as $V_1(\theta, \mu) = \phi_1 v_1(\theta, \mu) + \phi_2 v_2(\theta)$, where $v_1(\theta, \mu)$ and $v_2(\theta)$ are as defined in Section 2. Then the block Gibbs sampler satisfies (7) with*

$$b = \phi_1 \left[ c_1 + c_2 \sum_{i=1}^{K} m_i^{-1} + K\Delta^2 \right] + \phi_2 [c_2(K+1) + M\Delta^2].$$

PROOF. It suffices to show that

(14) $\quad E[V_1(\theta, \mu) | \lambda', \theta', \mu'] \leq \phi_1 \delta_1 v_1(\theta', \mu') + \left( \frac{\phi_1 \delta_4}{\phi_2} + \delta_3 \right) \phi_2 v_2(\theta') + b$

because

$$\phi_1 \delta_1 v_1(\theta', \mu') + \left( \frac{\phi_1 \delta_4}{\phi_2} + \delta_3 \right) \phi_2 v_2(\theta') + b$$
$$\leq \phi_1 \delta v_1(\theta', \mu') + \left( \frac{\phi_1 \delta_4}{\phi_2} + \delta \right) \phi_2 v_2(\theta') + b$$
$$\leq \gamma \phi_1 v_1(\theta', \mu') + \gamma \phi_2 v_2(\theta') + b$$
$$= \gamma V_1(\theta', \mu') + b.$$

In bounding the left-hand side of (14), we will use the following rule:

(15) $\quad E[V_1(\theta, \mu) | \lambda', \theta', \mu'] = E[V_1(\theta, \mu) | \theta', \mu'] = E\{E[V_1(\theta, \mu) | \lambda] | \theta', \mu'\},$

which follows from the form of the MTD for the block Gibbs sampler given in (5). We begin with some preliminary calculations. First, note that

$$E(\lambda_\theta^{-1}|\theta',\mu') = \frac{2b_1}{2a_1 + K - 2} + \frac{v_1(\theta',\mu')}{2a_1 + K - 2}$$
(16)
$$= c_1 + \delta_1 v_1(\theta',\mu')$$

and

$$E(\lambda_e^{-1}|\theta',\mu') = \frac{2b_2 + \text{SSE}}{2a_2 + M - 2} + \frac{v_2(\theta')}{2a_2 + M - 2}$$
(17)
$$= c_2 + \delta_2 v_2(\theta').$$

We now begin the main calculation. Using our rule, we have

$$E[v_1(\theta,\mu)|\theta',\mu'] = \sum_{i=1}^{K} E[(\theta_i - \mu)^2|\theta',\mu']$$

$$= E\left\{\sum_{i=1}^{K} E[(\theta_i - \mu)^2|\lambda]\Big|\theta',\mu'\right\}.$$

Using results from Appendix A, we have

$$E[(\theta_i - \mu)^2|\lambda]$$
$$= \text{Var}(\theta_i|\lambda) + \text{Var}(\mu|\lambda) - 2\text{Cov}[(\theta_i,\mu)|\lambda] + [E(\theta_i|\lambda) - E(\mu|\lambda)]^2$$
$$= \frac{1}{\lambda_\theta + m_i\lambda_e} + \frac{\lambda_\theta^2 + (\lambda_\theta + m_i\lambda_e)^2 - 2\lambda_\theta(\lambda_\theta + m_i\lambda_e)}{(s_0 + t)(\lambda_\theta + m_i\lambda_e)^2}$$
$$\quad + [E(\theta_i|\lambda) - E(\mu|\lambda)]^2$$
$$= \frac{1}{\lambda_\theta + m_i\lambda_e} + \frac{m_i^2\lambda_e^2}{(s_0 + t)(\lambda_\theta + m_i\lambda_e)^2} + [E(\theta_i|\lambda) - E(\mu|\lambda)]^2$$
$$\leq \frac{1}{m_i\lambda_e} + \frac{m_i\lambda_e}{t(\lambda_\theta + m_i\lambda_e)} + \Delta^2.$$

Hence,

(18) $$\sum_{i=1}^{K} E[(\theta_i - \mu)^2|\lambda] \leq \lambda_e^{-1}\sum_{i=1}^{K} m_i^{-1} + \lambda_\theta^{-1} + K\Delta^2.$$

Thus, by combining (16)–(18) we obtain

$$E[\phi_1 v_1(\theta,\mu)|\theta',\mu']$$
(19)
$$\leq \delta_1\phi_1 v_1(\theta',\mu') + \delta_4\phi_1 v_2(\theta') + \phi_1\left[c_1 + c_2\sum_{i=1}^{K} m_i^{-1} + K\Delta^2\right].$$

Now

$$E[v_2(\theta)|\theta',\mu'] = \sum_i m_i E[(\theta_i - \bar{y}_i)^2|\theta',\mu'] = E\bigg\{\sum_i m_i E[(\theta_i - \bar{y}_i)^2|\lambda]\Big|\theta',\mu'\bigg\}.$$

We can bound the innermost expectation as follows:

$$\begin{aligned}
E[(\theta_i - \bar{y}_i)^2|\lambda] &= \operatorname{Var}(\theta_i|\lambda) + [E(\theta_i|\lambda) - \bar{y}_i]^2 \\
&= \frac{1}{\lambda_\theta + m_i\lambda_e} + \frac{\lambda_\theta^2}{(s_0+t)(\lambda_\theta + m_i\lambda_e)^2} + [E(\theta_i|\lambda) - \bar{y}_i]^2 \\
&\leq \frac{1}{m_i\lambda_e} + \frac{\lambda_\theta}{t(\lambda_\theta + m_i\lambda_e)} + \Delta^2.
\end{aligned}$$

Hence

$$(20) \qquad \sum_{i=1}^{K} m_i E[(\theta_i - \bar{y}_i)^2|\lambda] \leq (K+1)\lambda_e^{-1} + M\Delta^2,$$

and so by combining (17) and (20) we obtain

$$(21) \qquad E[\phi_2 v_2(\theta)|\theta',\mu'] \leq \delta_3 \phi_2 v_2(\theta') + \phi_2[(K+1)c_2 + M\Delta^2].$$

Combining (19) and (21) yields (14). $\square$

REMARK 4.1. The upper bound on the total variation distance that is the conclusion of Theorem 3.1 involves the starting value of the Markov chain, $x_0$, *only* through $V(x_0)$. Moreover, given the way in which $V(x_0)$ enters the formula, it is clear that the optimal starting value, in terms of minimizing the upper bound, is the starting value that minimizes $V(x_0)$. This starting value is also optimal for the application of Theorem 3.2. In Appendix B we show that the value of $(\theta,\mu)$ that minimizes $V_1(\theta,\mu)$ has components

$$\hat{\theta}_i = \frac{\phi_1[\sum_{j=1}^{K}(m_j\bar{y}_j/(\phi_1+\phi_2 m_j))/\sum_{j=1}^{K}(m_j/(\phi_1+\phi_2 m_j))] + \phi_2 m_i \bar{y}_i}{\phi_1 + \phi_2 m_i}$$

and $\hat{\mu} = K^{-1}\sum_{i=1}^{K}\hat{\theta}_i$.

While the conclusion of Proposition 4.1 certainly holds when the data are balanced, it is possible to do better in this case. Specifically, the proof of Proposition 4.1 uses the general bounds on $[E(\theta_i|\lambda) - E(\mu|\lambda)]^2$ and $[E(\theta_i|\lambda) - \bar{y}_i]^2$ given in Appendix A. Much sharper bounds are possible by explicitly using the balancedness, and these lead to a better drift condition.

4.2. *Drift: balanced case.* Now assume that $m_i = m \geq 2$ for all $i = 1, \ldots, K$ and let $\delta_5 = K\delta_2 \in (0, 1)$.

PROPOSITION 4.2. *Fix $\gamma \in (\delta, 1)$ and let $\phi$ be a positive number such that $\phi\delta_5 + \delta < \gamma$. Define the drift function as $V_2(\theta, \mu) = \phi v_1(\theta, \mu) + m^{-1} v_2(\theta)$. Then the block Gibbs sampler satisfies (7) with*

$$b = \phi c_1 + [(\phi K + K + 1)/m]c_2 + \max\{\phi, 1\} \sum_{i=1}^{K} \max\{(\bar{y} - \bar{y}_i)^2, (m_0 - \bar{y}_i)^2\},$$

*where $\bar{y} := K^{-1} \sum_{i=1}^{K} \bar{y}_i$.*

PROOF. When the data are balanced,

$$t = \frac{M\lambda_\theta \lambda_e}{\lambda_\theta + m\lambda_e},$$

so that $E(\mu|\lambda) = (t\bar{y} + m_0 s_0)/(s_0 + t)$. Hence for all $i = 1, \ldots, K$ we have

$$[E(\theta_i|\lambda) - \bar{y}_i]^2 = \left[\frac{\lambda_\theta}{\lambda_\theta + m\lambda_e}\left(\frac{t\bar{y} + m_0 s_0}{s_0 + t}\right) + \frac{\lambda_e m \bar{y}_i}{\lambda_\theta + m\lambda_e} - \bar{y}_i\right]^2$$

$$= \left(\frac{\lambda_\theta}{\lambda_\theta + m\lambda_e}\right)^2 \left[\frac{t(\bar{y} - \bar{y}_i) + s_0(m_0 - \bar{y}_i)}{s_0 + t}\right]^2$$

$$\leq \left(\frac{\lambda_\theta}{\lambda_\theta + m\lambda_e}\right)^2 \frac{t(\bar{y} - \bar{y}_i)^2 + s_0(m_0 - \bar{y}_i)^2}{s_0 + t},$$

where the last inequality is Jensen's. A similar argument shows that, for all $i = 1, \ldots, K$,

$$[E(\theta_i|\lambda) - E(\mu|\lambda)]^2 \leq \left(\frac{m\lambda_e}{\lambda_\theta + m\lambda_e}\right)^2 \frac{t(\bar{y} - \bar{y}_i)^2 + s_0(m_0 - \bar{y}_i)^2}{s_0 + t}.$$

Therefore,

$$\phi[E(\theta_i|\lambda) - E(\mu|\lambda)]^2 + [E(\theta_i|\lambda) - \bar{y}_i]^2$$

$$\leq \max\{\phi, 1\}\left[\frac{t(\bar{y} - \bar{y}_i)^2 + s_0(m_0 - \bar{y}_i)^2}{s_0 + t}\right],$$

and hence

$$\sum_{i=1}^{K} \{\phi[E(\theta_i|\lambda) - E(\mu|\lambda)]^2 + [E(\theta_i|\lambda) - \bar{y}_i]^2\}$$

$$\leq \max\{\phi, 1\} \sum_{i=1}^{K} \max\{(\bar{y} - \bar{y}_i)^2, (m_0 - \bar{y}_i)^2\}.$$

To prove the result, it suffices to show that

(22) $\quad E[V_2(\theta,\mu)|\lambda',\theta',\mu'] \leq \phi\delta_1 v_1(\theta',\mu') + (\phi\delta_5 + \delta_3)m^{-1}v_2(\theta') + b$

since

$$\phi\delta_1 v_1(\theta',\mu') + (\phi\delta_5 + \delta_3)m^{-1}v_2(\theta') + b$$
$$\leq \phi\delta v_1(\theta',\mu') + (\phi\delta_5 + \delta)m^{-1}v_2(\theta') + b$$
$$\leq \gamma\phi v_1(\theta',\mu') + \gamma m^{-1}v_2(\theta') + b$$
$$= \gamma V_2(\theta',\mu') + b.$$

The remainder of the proof is nearly identical to the proof of Proposition 4.1 and is therefore left to the reader. □

REMARK 4.2. This result is stated (without proof) in Jones and Hobert [(2001), Appendix A] and the statement contains an error. Specifically, $b$ is stated incorrectly and should appear as above.

4.3. *Minorization.* We now use a technique based on Rosenthal's (1995a) Lemma 6b to establish a minorization condition of the form (8) on the set

$$S_B = \{(\theta,\mu) : V_1(\theta,\mu) \leq d\} = \{(\theta,\mu) : \phi_1 v_1(\theta,\mu) + \phi_2 v_2(\theta) \leq d\},$$

for any $d > 0$. Since $V_2$ of Proposition 4.2 is a special case of $V_1$, this minorization will also work for $V_2$. First note that $S_B$ is contained in $C_B := C_{B_1} \cap C_{B_2}$, where

$$C_{B_1} = \{(\theta,\mu) : v_1(\theta,\mu) < d/\phi_1\} \quad \text{and} \quad C_{B_2} = \{(\theta,\mu) : v_2(\theta) < d/\phi_2\}.$$

Hence, it suffices to establish a minorization condition that holds on $C_B$. We will accomplish this by finding an $\varepsilon > 0$ and a density $q(\lambda,\theta,\mu)$ on $\mathbb{R}_+^2 \times \mathbb{R}^K \times \mathbb{R}$ such that

$$k(\lambda,\theta,\mu|\lambda',\theta',\mu') \geq \varepsilon q(\lambda,\theta,\mu) \qquad \forall (\theta',\mu') \in C_B,$$

where $k(\lambda,\theta,\mu|\lambda',\theta',\mu')$ is the MTD for the block Gibbs sampler given in (5). We will require the following lemma, whose proof is given in Appendix C.

LEMMA 4.1. *Let* $\mathrm{Gamma}(\alpha,\beta;x)$ *denote the value of the* $\mathrm{Gamma}(\alpha,\beta)$ *density at the point* $x > 0$. *If* $\alpha > 1$, $b > 0$ *and* $c > 0$ *are fixed, then, as a function of* $x$,

$$\inf_{0<\beta<c} \mathrm{Gamma}(\alpha, b+\beta/2; x) = \begin{cases} \mathrm{Gamma}(\alpha, b; x), & \text{if } x < x^*, \\ \mathrm{Gamma}(\alpha, b+c/2; x), & \text{if } x > x^*, \end{cases}$$

*where*

$$x^* = \frac{2\alpha}{c}\log\left(1 + \frac{c}{2b}\right).$$

Here is the minorization condition.

PROPOSITION 4.3. *Let $q(\lambda, \theta, \mu)$ be a density on $\mathbb{R}_+^2 \times \mathbb{R}^K \times \mathbb{R}$ defined as*

$$q(\lambda, \theta, \mu) = \left[\frac{h_1(\lambda_\theta)}{\int_{\mathbb{R}_+} h_1(\lambda_\theta)\, d\lambda_\theta}\right]\left[\frac{h_2(\lambda_e)}{\int_{\mathbb{R}_+} h_2(\lambda_e)\, d\lambda_e}\right] f(\xi|\lambda, y),$$

*where*

$$h_1(\lambda_\theta) = \begin{cases} \operatorname{Gamma}\left(\dfrac{K}{2} + a_1, b_1; \lambda_\theta\right), & \lambda_\theta < \lambda_\theta^*, \\ \operatorname{Gamma}\left(\dfrac{K}{2} + a_1, \dfrac{d}{2\phi_1} + b_1; \lambda_\theta\right), & \lambda_\theta \geq \lambda_\theta^*, \end{cases}$$

*for*

$$\lambda_\theta^* = \frac{\phi_1(K + 2a_1)}{d} \log\left(1 + \frac{d}{2b_1\phi_1}\right)$$

*and*

$$h_2(\lambda_e) = \begin{cases} \operatorname{Gamma}\left(\dfrac{M}{2} + a_2, \dfrac{\mathrm{SSE}}{2} + b_2; \lambda_e\right), & \lambda_e < \lambda_e^*, \\ \operatorname{Gamma}\left(\dfrac{M}{2} + a_2, \dfrac{\phi_2 \mathrm{SSE} + d}{2\phi_2} + b_2; \lambda_e\right), & \lambda_e \geq \lambda_e^*, \end{cases}$$

*for*

$$\lambda_e^* = \frac{\phi_2(M + 2a_2)}{d} \log\left(1 + \frac{d}{\phi_2(2b_2 + \mathrm{SSE})}\right).$$

*Set $\varepsilon_B = [\int_{\mathbb{R}_+} h_1(\lambda_\theta)\, d\lambda_\theta][\int_{\mathbb{R}_+} h_2(\lambda_e)\, d\lambda_e]$. Then the Markov transition density for the block Gibbs sampler satisfies the following minorization condition:*

$$k(\lambda, \theta, \mu | \lambda', \theta', \mu') \geq \varepsilon_B q(\lambda, \theta, \mu) \qquad \forall\, (\theta', \mu') \in C_B.$$

PROOF. We use $\xi = (\theta, \mu)$ and $\xi' = (\theta', \mu')$ to simplify notation. If $\xi' \in C_B$, we have

$$\begin{aligned}
&f(\lambda_\theta|\xi', y) f(\lambda_e|\xi', y) f(\xi|\lambda, y) \\
&\quad \geq f(\xi|\lambda, y) \inf_{\xi \in C_B} [f(\lambda_\theta|\xi, y) f(\lambda_e|\xi, y)] \\
&\quad \geq f(\xi|\lambda, y) \left[\inf_{\xi \in C_B} f(\lambda_\theta|\xi, y)\right] \left[\inf_{\xi \in C_B} f(\lambda_e|\xi, y)\right] \\
&\quad \geq f(\xi|\lambda, y) \left[\inf_{\xi \in C_{B_1}} f(\lambda_\theta|\xi, y)\right] \left[\inf_{\xi \in C_{B_2}} f(\lambda_e|\xi, y)\right].
\end{aligned}$$

Thus we can take

$$q(\lambda, \theta, \mu) \propto f(\xi|\lambda, y) \left[\inf_{\xi \in C_{B_1}} f(\lambda_\theta|\xi, y)\right] \left[\inf_{\xi \in C_{B_2}} f(\lambda_e|\xi, y)\right].$$

Two applications of Lemma 4.1 yield the result. □

The drift and minorization conditions given in Propositions 4.1–4.3 can be used in conjunction with either Theorem 3.1 or 3.2 to get a formula giving an upper bound on the total variation distance to stationarity for the block Gibbs sampler. One such formula is stated explicitly at the start of Section 6.

**5. Drift and minorization for the Gibbs sampler.** In this section we develop drift and minorization conditions for the Gibbs sampler. We continue to assume that $m' = \min\{m_1, m_2, \ldots, m_K\} \geq 2$ and that $K \geq 3$. Let $m'' = \max\{m_1, m_2, \ldots, m_K\}$.

5.1. *Drift.* Recall that $\delta_1 = 1/(2a_1 + K - 2)$ and define

$$\delta_6 = \frac{K^2 + 2Ka_1}{2s_0 b_1 + K^2 + 2Ka_1} \quad \text{and} \quad \delta_7 = \frac{1}{2(a_1 - 1)}.$$

Clearly $\delta_6 \in (0, 1)$. It is straightforward to show that if $a_1 > 3/2$, then $\delta_7 \in (0, 1)$ and there exists $\rho_1 \in (0, 1)$ such that

$$\left(K + \frac{\delta_6}{\delta_7}\right)\delta_1 < \rho_1. \tag{23}$$

Define the function $v_3(\theta, \lambda) = \frac{K\lambda_\theta}{s_0 + K\lambda_\theta}(\bar{\theta} - \bar{y})^2$. Also, let $s^2 = \sum_{i=1}^{K}(\bar{y}_i - \bar{y})^2$. We will require the following lemma, whose proof is given in Appendix D.

LEMMA 5.1. *Let $a$ and $b$ be constants such that $5b > a \geq b > 0$. Then if $x$ and $y$ are positive,*

$$\left(\frac{ax}{ax + y}\right)^2 + \left(\frac{y}{bx + y}\right)^2 < 1. \tag{24}$$

Here is the drift condition.

PROPOSITION 5.1. *Assume that $a_1 > 3/2$ and let $\rho_1 \in (0, 1)$ satisfy (23). Assume also that $5m' > m''$. Fix $c_3 \in (0, \min\{b_1, b_2\})$ and fix $\gamma \in (\max\{\rho_1, \delta_6, \delta_7\}, 1)$. Define the drift function as*

$$V_3(\theta, \lambda) = e^{c_3\lambda_\theta} + e^{c_3\lambda_e} + \frac{\delta_7}{K\delta_1\lambda_\theta} + v_3(\theta, \lambda).$$

*Then the Gibbs sampler satisfies (7) with*

$$b = \left(\frac{b_1}{b_1 - c_3}\right)^{a_1+K/2} + \left(\frac{b_2}{b_2 - c_3}\right)^{a_2+N/2}$$

$$+ (\delta_6 + \delta_7)\left[\frac{1}{s_0} + (m_0 - \bar{y})^2 + \frac{s^2}{K}\right] + \frac{2b_1\delta_7}{K}.$$

PROOF. It suffices to show that

$$E[V_3(\theta, \lambda) | \mu', \theta', \lambda']$$

$$(25) \qquad \leq \frac{\delta_7[(K + \delta_6/\delta_7)\delta_1]}{K\delta_1\lambda'_\theta}$$

$$+ \left[\delta_7\left(\frac{m''\lambda'_e}{\lambda'_\theta + m''\lambda'_e}\right)^2 + \delta_6\left(\frac{\lambda'_\theta}{\lambda'_\theta + m'\lambda'_e}\right)^2\right] v_3(\theta', \lambda') + b,$$

because, using Lemma 5.1 and (23), we have

$$\frac{\delta_7}{K\delta_1\lambda'_\theta}\left[\left(K + \frac{\delta_6}{\delta_7}\right)\delta_1\right] + \left[\delta_7\left(\frac{m''\lambda'_e}{\lambda'_\theta + m''\lambda'_e}\right)^2 + \delta_6\left(\frac{\lambda'_\theta}{\lambda'_\theta + m'\lambda'_e}\right)^2\right] v_3(\theta', \lambda') + b$$

$$\leq \frac{\rho_1\delta_7}{K\delta_1\lambda'_\theta} + \max\{\delta_6, \delta_7\} v_3(\theta', \lambda') + b$$

$$\leq \gamma V_3(\theta', \lambda') + b.$$

Recall that we are considering Hobert and Geyer's (1998) updating scheme for the Gibbs sampler: $(\mu', \theta', \lambda') \to (\mu, \theta, \lambda)$. Establishing (25) requires the calculation of several expectations, and these will be calculated using the following rule:

$$E[V_3(\theta, \lambda) | \mu', \theta', \lambda'] = E[V_3(\theta, \lambda) | \theta', \lambda']$$
$$= E\{E\{E[V_3(\theta, \lambda) | \mu, \theta, \theta', \lambda'] | \mu, \theta', \lambda'\} | \theta', \lambda'\}$$
$$= E\{E\{E[V_3(\theta, \lambda) | \mu, \theta] | \mu, \lambda'\} | \theta', \lambda'\}.$$

We now establish (25). First, it is easy to show that

$$(26) \qquad E[e^{c_3\lambda_\theta} | \theta, \mu] \leq \left(\frac{b_1}{b_1 - c_3}\right)^{a_1+K/2} \quad \text{and}$$

$$E[e^{c_3\lambda_e} | \theta, \mu] \leq \left(\frac{b_2}{b_2 - c_3}\right)^{a_2+N/2}.$$

Now we evaluate $E[\frac{\delta_7}{K\delta_1\lambda_\theta} | \mu', \theta', \lambda']$. Note that

$$(27) \qquad E[\lambda_\theta^{-1} | \mu, \theta] = \delta_1\left[2b_1 + \sum_{i=1}^{K}(\theta_i - \mu)^2\right]$$

and

$$E[(\theta_i - \mu)^2|\mu, \lambda'] = \text{Var}(\theta_i|\mu, \lambda') + [E(\theta_i|\mu, \lambda') - \mu]^2$$

(28)
$$= \frac{1}{\lambda'_\theta + m_i\lambda'_e} + \left(\frac{m_i\lambda'_e}{\lambda'_\theta + m_i\lambda'_e}\right)^2 (\mu - \bar{y}_i)^2$$

$$\leq \frac{1}{\lambda'_\theta} + \left(\frac{m''\lambda'_e}{\lambda'_\theta + m''\lambda'_e}\right)^2 (\mu - \bar{y}_i)^2.$$

It follows that

(29) $$\sum_{i=1}^{K} E[(\theta_i - \mu)^2|\mu, \lambda'] \leq \frac{K}{\lambda'_\theta} + \left(\frac{m''\lambda'_e}{\lambda'_\theta + m''\lambda'_e}\right)^2 K(\mu - \bar{y})^2 + s^2.$$

Letting $\overline{\theta'} = K^{-1}\sum_i \theta'_i$, we have

$$E[(\mu - \bar{y})^2|\theta', \lambda'] = \text{Var}(\mu|\theta', \lambda') + [E(\mu|\theta', \lambda') - \bar{y}]^2$$

(30)
$$= \frac{1}{s_0 + K\lambda'_\theta} + \left[\frac{s_0}{s_0 + K\lambda'_\theta}(m_0 - \bar{y}) + \frac{K\lambda'_\theta}{s_0 + K\lambda'_\theta}(\overline{\theta'} - \bar{y})\right]^2$$

$$\leq \frac{1}{s_0 + K\lambda'_\theta} + \frac{s_0}{s_0 + K\lambda'_\theta}(m_0 - \bar{y})^2 + \frac{K\lambda'_\theta}{s_0 + K\lambda'_\theta}(\overline{\theta'} - \bar{y})^2$$

$$\leq \frac{1}{s_0} + (m_0 - \bar{y})^2 + v_3(\theta', \lambda'),$$

where the first inequality is Jensen's. On combining (27)–(30), we have

(31)
$$E\left[\frac{\delta_7}{K\delta_1\lambda_\theta}\Big|\mu', \theta', \lambda'\right] \leq \frac{\delta_7}{\lambda'_\theta} + \delta_7\left(\frac{m''\lambda'_e}{\lambda'_\theta + m''\lambda'_e}\right)^2 v_3(\theta', \lambda')$$

$$+ \delta_7\left[\frac{1}{s_0} + (m_0 - \bar{y})^2 + \frac{s^2}{K}\right] + \frac{2b_1\delta_7}{K}.$$

The last thing we need to evaluate is $E[v_3(\theta, \lambda)|\mu', \theta', \lambda']$. As in Hobert and Geyer (1998), Jensen's inequality yields

(32) $$E\left(\frac{K\lambda_\theta}{s_0 + K\lambda_\theta}\Big|\mu, \theta\right) \leq \frac{KE(\lambda_\theta|\mu, \theta)}{s_0 + KE(\lambda_\theta|\mu, \theta)} \leq \frac{K^2 + 2Ka_1}{2s_0b_1 + K^2 + 2Ka_1} = \delta_6.$$

These authors also note that the conditional independence of the $\theta_i$'s implies that

$$\bar{\theta}|\mu, \lambda \sim N\left(\frac{1}{K}\sum_{i=1}^{K}\frac{\lambda_\theta\mu + m_i\lambda_e\bar{y}_i}{\lambda_\theta + m_i\lambda_e}, \frac{1}{K^2}\sum_{i=1}^{K}\frac{1}{\lambda_\theta + m_i\lambda_e}\right),$$

from which it follows that

$$E[(\bar{\theta} - \bar{y})^2 | \mu, \lambda'] = \text{Var}(\bar{\theta} | \mu, \lambda') + [E(\bar{\theta} | \mu, \lambda') - \bar{y}]^2$$

$$= \frac{1}{K^2} \sum_{i=1}^{K} \frac{1}{\lambda'_\theta + m_i \lambda'_e} + \left[ \frac{1}{K} \sum_{i=1}^{K} \frac{\lambda'_\theta}{\lambda'_\theta + m_i \lambda'_e} (\mu - \bar{y}_i) \right]^2$$

(33)
$$\leq \frac{1}{K} \frac{1}{\lambda'_\theta + m' \lambda'_e} + \frac{1}{K} \sum_{i=1}^{K} \left( \frac{\lambda'_\theta}{\lambda'_\theta + m_i \lambda'_e} \right)^2 (\mu - \bar{y}_i)^2$$

$$\leq \frac{1}{K \lambda'_\theta} + \left( \frac{\lambda'_\theta}{\lambda'_\theta + m' \lambda'_e} \right)^2 (\mu - \bar{y})^2 + \frac{s^2}{K},$$

where, again, (part of) the first inequality is Jensen's. On combining (30), (32) and (33), we have

(34)
$$E[v_3(\theta, \lambda) | \mu', \theta', \lambda']$$
$$\leq \frac{\delta_6}{K \lambda'_\theta} + \delta_6 \left( \frac{\lambda'_\theta}{\lambda'_\theta + m' \lambda'_e} \right)^2 v_3(\theta', \lambda') + \delta_6 \left[ \frac{1}{s_0} + (m_0 - \bar{y})^2 + \frac{s^2}{K} \right].$$

Combining (26), (31) and (34) yields (25). $\square$

REMARK 5.1. Note that our drift condition for the block Gibbs sampler (Proposition 4.1) holds for all hyperparameter configurations (corresponding to proper priors) and nearly all values of $m'$ and $m''$. In contrast, it is assumed in Proposition 5.1 that $a_1 > 3/2$ and that $5m' > m''$. On the other hand, Hobert and Geyer's (1998) drift condition for the Gibbs sampler involves even more restrictive assumptions about $a_1$ and the relationship between $m'$ and $m''$. Specifically, Hobert and Geyer (1998) assume that $a_1 \geq (3K - 2)/(2K - 2)$ and that $m' > (\sqrt{5} - 2)m''$. Note that $(3K - 2)/(2K - 2) > 3/2$ for all $K \geq 2$ and that $5 > (\sqrt{5} - 2)^{-1} \approx 4.23$.

REMARK 5.2. In this case the optimal starting value minimizes

$$V_3(\theta, \lambda) = e^{c_3 \lambda_\theta} + e^{c_3 \lambda_e} + \frac{\delta_7}{K \delta_1 \lambda_\theta} + \frac{K \lambda_\theta}{s_0 + K \lambda_\theta} (\bar{\theta} - \bar{y})^2.$$

The last term will vanish as long as the $\theta_i$'s are such that $\bar{\theta} = \bar{y}$. The optimal starting value for $\lambda_\theta$ is the minimizer of the function $e^{c_3 \lambda_\theta} + \delta_7 / (K \delta_1 \lambda_\theta)$. This cannot be computed in closed form, but is easily found numerically. Finally, since $\lambda_e = 0$ is not appropriate, we simply start $\lambda_e$ at a small positive number.

5.2. *Minorization.* Fix $d > 0$ and define $S_G = \{(\theta, \lambda) : V_3(\theta, \lambda) \leq d\}$. Similar to our previous work with the block Gibbs sampler, our goal will be to find a density $q(\mu, \theta, \lambda)$ on $\mathbb{R} \times \mathbb{R}^K \times \mathbb{R}^2_+$ and an $\varepsilon > 0$ such that

$$k(\mu, \theta, \lambda | \mu', \lambda', \theta') \geq \varepsilon q(\mu, \theta, \lambda) \qquad \forall (\theta', \lambda') \in S_G.$$

As before, we will actually establish the minorization on a superset of $S_G$ with which it is more convenient to work. Let $c_4 = \delta_7/(K\delta_1 d)$ and put $c_l$ and $c_u$ equal to $\bar{y} - \sqrt{(m_0 - \bar{y})^2 + d}$ and $\bar{y} + \sqrt{(m_0 - \bar{y})^2 + d}$, respectively. We show in Appendix E that $S_G \subset C_G = C_{G_1} \cap C_{G_2} \cap C_{G_3}$, where

$$C_{G_1} = \left\{(\theta, \lambda) : c_4 \leq \lambda_\theta \leq \frac{\log d}{c_3}\right\}, \qquad C_{G_2} = \left\{(\theta, \lambda) : 0 < \lambda_e \leq \frac{\log d}{c_3}\right\},$$

$$C_{G_3} = \left\{(\theta, \lambda) : c_l \leq \frac{s_0 m_0 + K\lambda_\theta \bar{\theta}}{s_0 + K\lambda_\theta} \leq c_u\right\}.$$

Also, $C_{G_1} \cap C_{G_2}$ is nonempty as long as $d\log d > (c_3 \delta_7)/(K\delta_1)$. We will require the following obvious lemma.

LEMMA 5.2. *Let* $\mathrm{N}(\tau, \sigma^2; x)$ *denote the value of the* $\mathrm{N}(\tau, \sigma^2)$ *density at the point* $x$. *If* $a \leq b$, *then, as a function of* $x$,

$$\inf_{a \leq \tau \leq b} \mathrm{N}(\tau, \sigma^2; x) = \begin{cases} \mathrm{N}(b, \sigma^2; x), & \text{if } x \leq (a+b)/2, \\ \mathrm{N}(a, \sigma^2; x), & \text{if } x > (a+b)/2. \end{cases}$$

Here is the minorization condition.

PROPOSITION 5.2. *Let* $q(\mu, \theta, \lambda)$ *be a density on* $\mathbb{R} \times \mathbb{R}^K \times \mathbb{R}_+^2$ *defined as follows:*

$$q(\mu, \theta, \lambda) = \left[\frac{g_1(\mu, \theta) g_2(\mu)}{\int_\mathbb{R} \int_{\mathbb{R}^K} g_1(\mu, \theta) g_2(\mu) \, d\theta \, d\mu}\right] f(\lambda | \mu, \theta, y),$$

*where*

$$g_1(\mu, \theta) = \left(\frac{c_4}{2\pi}\right)^{K/2} \exp\left\{-\frac{\log d}{2c_3} \sum_{i=1}^{K} [(\theta_i - \mu)^2 + m_i(\theta_i - \bar{y}_i)^2]\right\}$$

*and*

$$g_2(\mu) = \begin{cases} \mathrm{N}(c_u, [s_0 + K\log(d)/c_3]^{-1}; \mu), & \mu \leq \bar{y}, \\ \mathrm{N}(c_l, [s_0 + K\log(d)/c_3]^{-1}; \mu), & \mu > \bar{y}. \end{cases}$$

*Set*

$$\varepsilon_G = \left[\frac{s_0 + Kc_4}{s_0 + K\log(d)/c_3}\right]^{1/2} \left[\int_\mathbb{R} \int_{\mathbb{R}^K} g_1(\mu, \theta) g_2(\mu) \, d\theta \, d\mu\right].$$

*Then the Markov transition density for the Gibbs sampler satisfies the minorization condition*

$$k(\mu, \theta, \lambda | \mu', \theta', \lambda') \geq \varepsilon_G q(\mu, \theta, \lambda) \qquad \forall (\theta', \lambda') \in C_G.$$

PROOF. Recall that $k(\mu, \theta, \lambda | \mu', \theta', \lambda') = f(\mu | \theta', \lambda', y) f(\theta | \mu, \lambda', y) f(\lambda | \mu, \theta, y)$.
For $(\theta', \lambda') \in C_G$, we have

$$f(\mu|\theta', \lambda', y) f(\theta|\mu, \lambda', y)$$
$$\geq \inf_{(\theta', \lambda') \in C_G} f(\mu|\theta', \lambda', y) f(\theta|\mu, \lambda', y)$$
$$\geq \left[\inf_{(\theta', \lambda') \in C_G} f(\mu|\theta', \lambda', y)\right] \left[\inf_{(\theta', \lambda') \in C_G} f(\theta|\mu, \lambda', y)\right]$$
$$\geq \left[\inf_{(\theta', \lambda') \in C_G} f(\mu|\theta', \lambda', y)\right] \left[\inf_{\lambda' \in C_{G_1} \cap C_{G_2}} f(\theta|\mu, \lambda', y)\right].$$

Using the fact that the $\theta_i$'s are conditionally independent, we have

$$\inf_{\lambda' \in C_{G_1} \cap C_{G_2}} f(\theta|\mu, \lambda', y)$$
$$= \inf_{\lambda' \in C_{G_1} \cap C_{G_2}} \prod_{i=1}^K f(\theta_i|\mu, \lambda', y) \geq \prod_{i=1}^K \inf_{\lambda' \in C_{G_1} \cap C_{G_2}} f(\theta_i|\mu, \lambda', y).$$

Now, using Jensen's inequality again, we have

$$f(\theta_i|\mu, \lambda', y)$$
$$= \sqrt{\frac{\lambda'_\theta + m_i \lambda'_e}{2\pi}} \exp\left\{-\frac{\lambda'_\theta + m_i \lambda'_e}{2}\left(\theta_i - \frac{\lambda'_\theta \mu + m_i \lambda'_e \bar{y}_i}{\lambda'_\theta + m_i \lambda'_e}\right)^2\right\}$$
$$= \sqrt{\frac{\lambda'_\theta + m_i \lambda'_e}{2\pi}}$$
$$\times \exp\left\{-\frac{\lambda'_\theta + m_i \lambda'_e}{2}\left[\frac{\lambda'_\theta}{\lambda'_\theta + m_i \lambda'_e}(\theta_i - \mu) + \frac{m_i \lambda'_e}{\lambda'_\theta + m_i \lambda'_e}(\theta_i - \bar{y}_i)\right]^2\right\}$$
$$\geq \sqrt{\frac{\lambda'_\theta + m_i \lambda'_e}{2\pi}}$$
$$\times \exp\left\{-\frac{\lambda'_\theta + m_i \lambda'_e}{2}\left[\frac{\lambda'_\theta}{\lambda'_\theta + m_i \lambda'_e}(\theta_i - \mu)^2 + \frac{m_i \lambda'_e}{\lambda'_\theta + m_i \lambda'_e}(\theta_i - \bar{y}_i)^2\right]\right\}$$
$$= \sqrt{\frac{\lambda'_\theta + m_i \lambda'_e}{2\pi}} \exp\left\{-\frac{1}{2}[\lambda'_\theta(\theta_i - \mu)^2 + m_i \lambda'_e(\theta_i - \bar{y}_i)^2]\right\}.$$

Hence,

$$\inf_{\lambda' \in C_{G_1} \cap C_{G_2}} f(\theta|\mu, \lambda', y)$$
$$\geq \left(\frac{c_4}{2\pi}\right)^{K/2} \exp\left\{-\frac{\log d}{2c_3} \sum_{i=1}^K [(\theta_i - \mu)^2 + m_i(\theta_i - \bar{y}_i)^2]\right\}$$
$$= g_1(\mu, \theta).$$

Now, if $(\theta', \lambda') \in C_G$, then $c_4 \leq \lambda'_\theta \leq \frac{\log d}{c_3}$ and hence

$$f(\mu|\theta', \lambda', y) = \sqrt{\frac{s_0 + K\lambda'_\theta}{2\pi}} \exp\left\{-\frac{s_0 + K\lambda'_\theta}{2}\left(\mu - \frac{s_0 m_0 + K\lambda'_\theta \overline{\theta'}}{s_0 + K\lambda'_\theta}\right)^2\right\}$$

$$\geq \sqrt{\frac{s_0 + Kc_4}{s_0 + (K\log d)/c_3}} \sqrt{\frac{s_0 + (K\log d)/c_3}{2\pi}}$$

$$\times \exp\left\{-\frac{s_0 + (K\log d)/c_3}{2}\left(\mu - \frac{s_0 m_0 + K\lambda'_\theta \overline{\theta'}}{s_0 + K\lambda'_\theta}\right)^2\right\}.$$

Thus,

$$\inf_{(\theta', \lambda') \in C_G} f(\mu|\theta', \lambda', y)$$

$$\geq \sqrt{\frac{s_0 + Kc_4}{s_0 + (K\log d)/c_3}} \inf_{(\theta', \lambda') \in C_G} \sqrt{\frac{s_0 + (K\log d)/c_3}{2\pi}}$$

$$\times \exp\left\{-\frac{s_0 + (K\log d)/c_3}{2}\left(\mu - \frac{s_0 m_0 + K\lambda'_\theta \overline{\theta'}}{s_0 + K\lambda'_\theta}\right)^2\right\}$$

$$\geq \sqrt{\frac{s_0 + Kc_4}{s_0 + (K\log d)/c_3}} \inf_{(\theta', \lambda') \in C_{G_3}} \sqrt{\frac{s_0 + (K\log d)/c_3}{2\pi}}$$

$$\times \exp\left\{-\frac{s_0 + (K\log d)/c_3}{2}\left(\mu - \frac{s_0 m_0 + K\lambda'_\theta \overline{\theta'}}{s_0 + K\lambda'_\theta}\right)^2\right\}$$

$$\geq g_2(\mu)\sqrt{\frac{s_0 + Kc_4}{s_0 + (K\log d)/c_3}},$$

where the last inequality is an application of Lemma 5.2. □

REMARK 5.3. In Appendix F we give a closed form expression for $\varepsilon_G$ involving the standard normal cumulative distribution function.

**6. A numerical example.** Consider a balanced data situation and let $\pi_h(\cdot)$ denote the probability measure corresponding to the posterior density in (2). Let $P^n((\lambda_0, \xi_0), \cdot)$ denote the $n$-step Markov transition kernel for the block Gibbs sampler started at $(\lambda_0, \xi_0)$. [Equation (5) shows that a starting value for $\lambda_0$ is actually not required.] We now write down an explicit upper bound for

$$\|P^n((\lambda_0, \xi_0), \cdot) - \pi_h(\cdot)\|,$$

based on Theorem 3.1 and Propositions 4.2 and 4.3. Although it has been suppressed in the notation, both $\pi_h$ and $P^n$ depend heavily on the six hyperparameters, $a_1$, $b_1$, $a_2$, $b_2$, $s_0$ and $m_0$. Our upper bound holds for all hyperparameter configurations such that $a_1, b_1, a_2, b_2, s_0$ are positive, that is, all hyperparameter configurations such that the priors on $\lambda_\theta$, $\lambda_e$ and $\mu$ are proper. Due to its generality, the bound is complicated to state. First, recall that $\text{SSE} = \sum_{i,j}(y_{ij} - \bar{y}_i)^2$, where $\bar{y}_i = m^{-1}\sum_{j=1}^m y_{ij}$. Recall further that

$$\delta_1 = \frac{1}{2a_1 + K - 2}, \qquad \delta_2 = \frac{1}{2a_2 + M - 2},$$

$\delta_3 = (K+1)\delta_2$, $\delta_5 = K\delta_2$, $\delta = \max\{\delta_1, \delta_3\}$, $c_1 = 2b_1\delta_1$ and $c_2 = (2b_2 + \text{SSE})\delta_2$. Note that all of these quantities depend only on the data and the hyperparameters.

Now choose $\gamma \in (\delta, 1)$ and $\phi > 0$ such that $\phi\delta_5 + \delta < \gamma$. Also, let

$$b = \phi c_1 + [(\phi K + K + 1)/m]c_2 + \max\{\phi, 1\}\sum_{i=1}^K \max\{(\bar{y} - \bar{y}_i)^2, (m_0 - \bar{y}_i)^2\},$$

and choose $d_R > 2b/(1 - \gamma)$. Finally, let

$$\varepsilon_B = \left[\int_{\mathbb{R}_+} h_1(\lambda_\theta) \, d\lambda_\theta\right]\left[\int_{\mathbb{R}_+} h_2(\lambda_e) \, d\lambda_e\right],$$

where

$$h_1(\lambda_\theta) = \begin{cases} \text{Gamma}\left(\dfrac{K}{2} + a_1, b_1; \lambda_\theta\right), & \lambda_\theta < \lambda_\theta^*, \\ \text{Gamma}\left(\dfrac{K}{2} + a_1, \dfrac{d_R}{2\phi} + b_1; \lambda_\theta\right), & \lambda_\theta \geq \lambda_\theta^*, \end{cases}$$

for

$$\lambda_\theta^* = \frac{\phi(K + 2a_1)}{d_R}\log\left(1 + \frac{d_R}{2b_1\phi}\right)$$

and

$$h_2(\lambda_e) = \begin{cases} \text{Gamma}\left(\dfrac{M}{2} + a_2, \dfrac{\text{SSE}}{2} + b_2; \lambda_e\right), & \lambda_e < \lambda_e^*, \\ \text{Gamma}\left(\dfrac{M}{2} + a_2, \dfrac{\text{SSE} + md_R}{2} + b_2; \lambda_e\right), & \lambda_e \geq \lambda_e^*, \end{cases}$$

for

$$\lambda_e^* = \frac{(M + 2a_2)}{md_R}\log\left(1 + \frac{md_R}{2b_2 + \text{SSE}}\right).$$

TABLE 1
*Simulated data*

| Cell | 1 | 2 | 3 | 4 | 5 |
|---|---|---|---|---|---|
| $\bar{y}_i$ | $-0.80247$ | $-1.0014$ | $-0.69090$ | $-1.1413$ | $-1.0125$ |

$$M = mK = 50$$
$$\bar{y} = M^{-1} \sum_{i=1}^{5} \sum_{j=1}^{10} y_{ij} = -0.92973$$
$$\text{SSE} = \sum_{i=1}^{5} \sum_{j=1}^{10} (y_{ij} - \bar{y}_i)^2 = 32.990$$

Note that $\varepsilon_\text{B}$ cannot be calculated in closed form, but can be evaluated numerically with four calls to a routine that evaluates the incomplete gamma function. Recall from the statement of Theorem 3.1 that

$$\alpha = \frac{1 + d_\text{R}}{1 + 2b + \gamma d_\text{R}} \quad \text{and} \quad U = 1 + 2(\gamma d_\text{R} + b).$$

Here is the bound. For any $0 < r < 1$ and any $n \in \{1, 2, 3, \ldots\}$,

$$\|P^n((\lambda_0, \xi_0), \cdot) - \pi_h(\cdot)\|$$
$$\leq (1 - \varepsilon_\text{B})^{rn} + \left(\frac{U^r}{\alpha^{1-r}}\right)^n \left(1 + \frac{b}{1-\gamma} + \phi v_1(\theta_0, \mu_0) + m^{-1} v_2(\theta_0)\right).$$

Using the optimal starting values from Remark 4.1, this becomes

$$\begin{aligned}(35) \quad &\|P^n((\lambda_0, \xi_0^\text{opt}), \cdot) - \pi_h(\cdot)\| \\ &\leq (1 - \varepsilon_\text{B})^{rn} + \left(\frac{U^r}{\alpha^{1-r}}\right)^n \left(1 + \frac{b}{1-\gamma} + \frac{\phi}{1+\phi} \sum_{i=1}^{K} (\bar{y}_i - \bar{y})^2\right).\end{aligned}$$

Explicit upper bounds can also be written for the block Gibbs sampler in the unbalanced case and for the Gibbs sampler. These are similar and are left to the reader. It is interesting to note that because our drift and minorization conditions for the block Gibbs sampler are free of $s_0$, so too is the bound in (35).

To evaluate (35), the user must provide values for $\gamma$, $\phi$, $d_\text{R}$ and $r$. In our experience, small changes in these quantities can lead to dramatically different results. Unfortunately, the right-hand side of (35) is a very complicated function of $\gamma$, $\phi$, $d_\text{R}$ and $r$. Hence, it would be quite difficult to find "optimal" values. In our applications of (35), we simply define reasonable ranges for these four quantities and then perform a grid search to find the configuration that leads to the smallest upper bound. We now provide an example of the use of (35) and of the analogous bound based on Theorem 3.2.

The data in Table 1 were simulated according to the model defined in Section 2 with $K = 5$, $m = 10$, $a_1 = 2.5$, $a_2 = b_1 = b_2 = 1$, $m_0 = 0$ and $s_0 = 1$.

We now pretend that the origin of the data is unknown and consider using the block Gibbs sampler to make approximate draws from four different intractable posterior distributions corresponding to the four hyperparameter settings listed in Table 2. The first setting in Table 2 is the "correct" prior in that it is exactly the setting under which the data were simulated. As one moves from setting 2 to setting 4, the prior variances on $\lambda_\theta$ and $\lambda_e$ become larger; that is, the priors become more "diffuse." For reasons discussed below $m_0$ is set equal to $\bar{y}$ in settings 2–4.

For each of the hyperparameter settings in Table 2 we used (35) as well as the analogous bound based on Theorem 3.2 to find an $n^*$ such that

$$\|P^{n^*}((\lambda_0, \xi_0^{\text{opt}}), \cdot) - \pi_h(\cdot)\| \leq 0.01.$$

The results are given in Tables 3 and 4. For example, consider hyperparameter setting 2. Theorem 3.1 yields

$$\|P^{3415}((\lambda_0, \xi_0^{\text{opt}}), \cdot) - \pi_h(\cdot)\| \leq 0.00999,$$

while Theorem 3.2 yields

$$\|P^{6563}((\lambda_0, \xi_0^{\text{opt}}), \cdot) - \pi_h(\cdot)\| \leq 0.00999.$$

While examining the $n^*$'s in Tables 3 and 4, keep in mind that it takes about 1.5 minutes to run one million iterations of the block Gibbs sampler on a standard PC. Thus, even the larger $n^*$'s are *feasible*.

Note that the results based on Theorem 3.1 are better across the board than those based on Theorem 3.2. We suspect that our use of Lemma 3.1 in the application of Theorem 3.2 has somewhat (artificially) inflated the $n^*$'s in Table 4.

A comparison of the $n^*$'s for hyperparameter settings 1 and 2 (in either table) shows that our bound is extremely sensitive to the distance between $m_0$ and $\bar{y}$. This is due to the fact that $\varepsilon_B$ decreases rapidly as $b$ increases and $b$ contains the term $\sum_{i=1}^{K} \max\{(\bar{y} - \bar{y}_i)^2, (m_0 - \bar{y}_i)^2\}$, which is minimized when $m_0 = \bar{y}$. While there may actually be some difference in the convergence rates of the two Markov chains corresponding to settings 1 and 2, it seems unlikely that the difference is as large as these numbers suggest. (Remember, these are only sufficient burn-ins.) It is probably the case that our results simply produce a better bound under setting 2 than they do under setting 1. This issue is discussed further in Section 7.

TABLE 2
*Four different prior specifications*

| Hyperparameter setting | $a_1$ | $b_1$ | $a_2$ | $b_2$ | $m_0$ |
|---|---|---|---|---|---|
| 1 | 2.5 | 1 | 1 | 1 | 0 |
| 2 | 2.5 | 1 | 1 | 1 | $\bar{y}$ |
| 3 | 0.1 | 0.1 | 0.1 | 0.1 | $\bar{y}$ |
| 4 | 0.01 | 0.01 | 0.01 | 0.01 | $\bar{y}$ |

TABLE 3
*Total variation bounds for the block Gibbs sampler via Theorem 3.1*

| Hyperparameter setting | $\gamma$ | $\phi$ | $d_{\mathrm{R}}$ | $r$ | $\varepsilon_{\mathrm{B}}$ | $n^*$ | Bound |
|---|---|---|---|---|---|---|---|
| 1 | 0.2596 | 0.9423 | 15.997 | 0.0188 | $3.1 \times 10^{-7}$ | $7.94 \times 10^8$ | 0.00999 |
| 2 | 0.2596 | 0.5385 | 3.0079 | 0.0789 | 0.0171 | $3.415 \times 10^3$ | 0.00999 |
| 3 | 0.4183 | 0.3059 | 2.8351 | 0.0512 | $6.8 \times 10^{-4}$ | $1.315 \times 10^5$ | 0.00999 |
| 4 | 0.4340 | 0.2965 | 2.8039 | 0.0483 | $8.1 \times 10^{-6}$ | $1.1796 \times 10^7$ | 0.00999 |

Another noteworthy feature of Tables 3 and 4 is that $n^*$ increases as the priors become more "diffuse." Figure 1 contains two plots describing the relationship between the prior variances on $\lambda_\theta$ and $\lambda_e$ and $n^*$. [The $n^*$'s in this plot were calculated using (35).] Note that $n^*$ increases quite rapidly with the prior variance on $\lambda_\theta$. While it is tempting to conclude that the chains associated with "diffuse" priors are relatively slow to converge, we cannot be sure that this is the case because, again, these are only sufficient burn-ins. However, our findings are entirely consistent with the work of Natarajan and McCulloch (1998), whose empirical results suggest that the mixing rate of the Gibbs sampler (for a probit–normal hierarchical model) becomes much slower as the priors become more diffuse.

**7. Discussion.** The quality of the upper bounds produced using Theorems 3.1 and 3.2 depends not only on the sharpness of the inequalities used to prove the theorems themselves, but also on the quality of the drift and minorization conditions used in the particular application. Consequently, it is possible, and perhaps even likely, that the chains we have analyzed actually get within 0.01 of stationarity much sooner than the $n^*$'s in Tables 3 and 4 would suggest. For example, we *know* from Table 3 that a sufficient burn-in for hyperparameter setting 2 is 3415. Thus, the value 6563 from Table 4 is too large by *at least* a factor of 1.9. The question then becomes how conservative are the results based on Rosenthal's theorem? As we now explain, this question was addressed by van Dyk and Meng (2001) in a different context.

TABLE 4
*Total variation bounds for the block Gibbs sampler via Theorem 3.2*

| Hyperparameter setting | $\rho$ | $\phi$ | $d_{\mathrm{RT}}$ | $\varepsilon_{\mathrm{B}}$ | $n^*$ | Bound |
|---|---|---|---|---|---|---|
| 1 | 0.615 | 0.84 | 15.213 | $4.1 \times 10^{-7}$ | $1.8835 \times 10^9$ | 0.00999 |
| 2 | 0.5975 | 0.49 | 2.6564 | 0.0234 | $6.563 \times 10^3$ | 0.00999 |
| 3 | 0.7113 | 0.3181 | 2.8492 | $7.2 \times 10^{-4}$ | $3.3915 \times 10^5$ | 0.00999 |
| 4 | 0.7191 | 0.3084 | 2.8154 | $8.6 \times 10^{-6}$ | $2.966 \times 10^7$ | 0.00999 |

Hobert (2001) used Theorem 3.1 to calculate a sufficient burn-in for a Markov chain Monte Carlo (MCMC) algorithm developed in Meng and van Dyk (1999). In the Rejoinder of van Dyk and Meng (2001) an empirical estimator of the total variation distance to stationarity was developed and used to demonstrate that Hobert's upper bound is probably extremely conser-

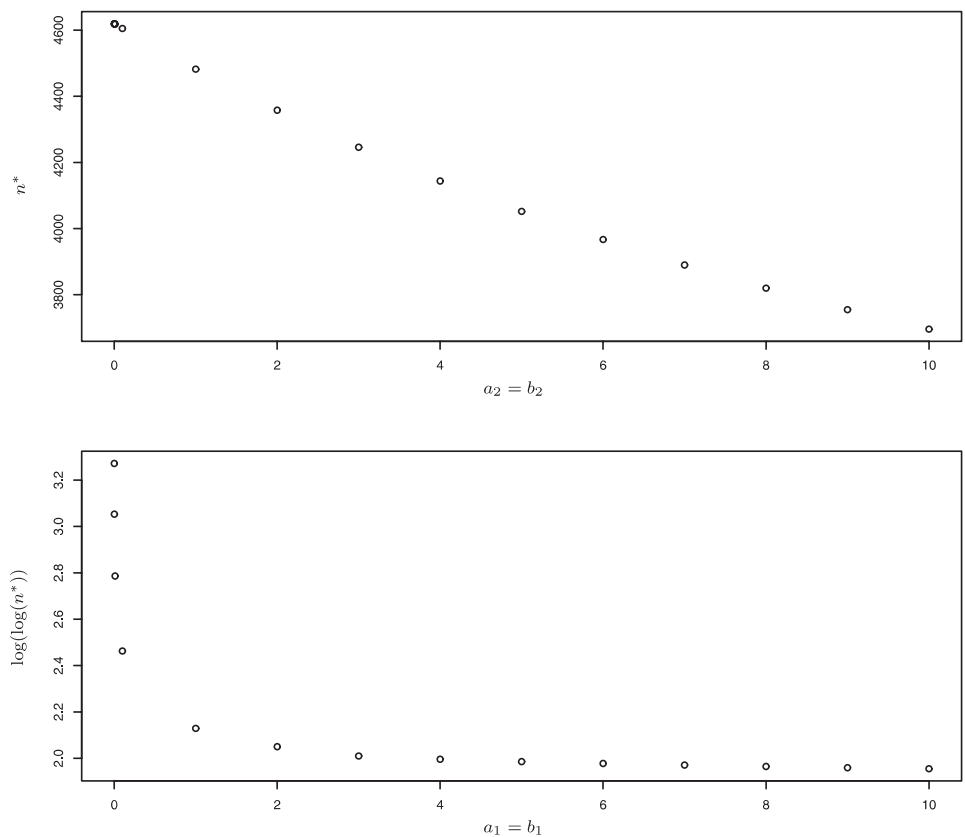

FIG. 1. *These two plots show how the "diffuseness" of the priors on $\lambda_\theta$ and $\lambda_e$ affects $n^*$. The top plot shows $n^*$ against $a_2 = b_2$ where the hyperparameters associated with $\lambda_\theta$ are held constant at $a_1 = b_1 = 1$. When $a_2 = b_2$, the prior variance of $\lambda_e$ is $1/b_2$ and the prior mean is constant at 1. The bottom plot shows $\log(\log(n^*))$ against $a_1 = b_1$ where the hyperparameters associated with $\lambda_e$ are held constant at $a_2 = b_2 = 1$. When $a_1 = b_1$, the prior variance of $\lambda_\theta$ is $1/b_1$ and the prior mean is constant at 1. In all cases $m_0$ was set equal to $\bar{y}$.*

TABLE 5
*Simulated data*

| Cell | 1 | 2 | 3 |
|------|---|---|---|
| $\bar{y}_i$ | $-0.54816$ | $0.92516$ | $-0.19924$ |

$$M_T = mK = 12$$

$$\bar{y} = M_T^{-1} \sum_{i=1}^{3} \sum_{j=1}^{4} y_{ij} = 0.059253$$

$$\text{SSE} = \sum_{i=1}^{3} \sum_{j=1}^{4} (y_{ij} - \bar{y}_i)^2 = 20.285$$

vative. Indeed, Hobert's sufficient burn-in was $n^* = 335$ while van Dyk and Meng's simulation results suggested that a burn-in of 2 is sufficient. We have experimented with van Dyk and Meng's empirical techniques in our situation and have come to similar conclusions. It would be interesting to use a Markov chain whose convergence behavior is known exactly to study how the sharpness of the bounds produced by Theorems 3.1 and 3.2 changes when different drift and minorization conditions are used.

In situations where it is possible to rigorously analyze two different MCMC algorithms for the same family of intractable posteriors, it is tempting to compare the algorithms using sufficient burn-in. However, we do not believe that this is an entirely fair method of comparison. Consider using our results in this way to compare Gibbs and block Gibbs. As we mentioned above, our Gibbs sampler is more difficult to analyze than our block Gibbs sampler. This probably results in relatively lower quality drift and minorization conditions for the Gibbs sampler. Indeed, using Propositions 5.1 and 5.2 in conjunction with Theorem 3.1 almost always yields extremely large $n^*$'s. Specifically, unless the priors are extremely "informative," it is difficult to find a hyperparameter configuration under which $\varepsilon_G$ is not effectively 0. Here is a comparison.

The data in Table 5 were simulated according to the model defined in Section 2 with $K = 3$, $m = 4$, $a_1 = a_2 = b_1 = b_2 = 2$, $s_0 = 1$ and $m_0 = 0$. We use the informative hyperparameter setting: $a_1 = 5$, $a_2 = 2$, $b_1 = 20$, $b_2 = 20$, $m_0 = 0$ and $s_0 = 4$. For the block Gibbs sampler (35) yields

$$\|P^{16631}((\lambda_0, \xi_0^{\text{opt}}), \cdot) - \pi_h(\cdot)\| \leq 0.00999.$$

For the Gibbs sampler Propositions 5.1 and 5.2 in conjunction with Theorem 3.1 yield

$$\|P_G^{4.826 \times 10^{19}}((\mu_0, \theta_0^{\text{opt}}, \lambda_0^{\text{opt}}), \cdot) - \pi_h(\cdot)\| \leq 0.00999.$$

As starting values for the Gibbs sampler we used $(\theta_0^{\text{opt}}, \lambda_0^{\text{opt}}) = (\bar{y}, \bar{y}, \bar{y}, 10^{-6}, 0.2839)$ (see Remark 5.2). The constants used to construct these bounds are given in Table 6.

While it is probably the case that block Gibbs converges faster than Gibbs, it is unlikely that the true difference is anywhere near as large as these numbers suggest. Thus, if we use these results to compare Gibbs and block Gibbs, the former will be penalized by the fact that it is simply more analytically cumbersome.

## APPENDIX A

**A.1. The elements of $\xi^*$ and $V$.** Hobert and Geyer [(1998), page 418] show that $\xi|\lambda, y \sim N(\xi^*, V)$ and give the specific forms of $\xi^* = \xi^*(\lambda, y)$ and $V = V(\lambda, y)$. We restate their results here. First we let

$$t = \sum_{i=1}^{K} \frac{m_i \lambda_\theta \lambda_e}{\lambda_\theta + m_i \lambda_e},$$

then

$$\text{Var}(\theta_i|\lambda) = \frac{1}{\lambda_\theta + m_i \lambda_e}\left[1 + \frac{\lambda_\theta^2}{(\lambda_\theta + m_i \lambda_e)(s_0 + t)}\right],$$

$$\text{Cov}(\theta_i, \theta_j|\lambda) = \frac{\lambda_\theta^2}{(\lambda_\theta + m_i \lambda_e)(\lambda_\theta + m_j \lambda_e)(s_0 + t)},$$

$$\text{Cov}(\theta_i, \mu|\lambda) = \frac{\lambda_\theta}{(\lambda_\theta + m_i \lambda_e)(s_0 + t)},$$

$$\text{Var}(\mu|\lambda) = \frac{1}{s_0 + t}.$$

Finally,

$$E(\mu|\lambda) = \frac{1}{s_0 + t}\left[\sum_{i=1}^{K} \frac{m_i \lambda_\theta \lambda_e \bar{y}_i}{\lambda_\theta + m_i \lambda_e} + m_0 s_0\right]$$

and

$$E(\theta_i|\lambda) = \frac{\lambda_\theta}{\lambda_\theta + m_i \lambda_e}\left[\frac{1}{s_0 + t}\left[\sum_{j=1}^{K} \frac{m_j \lambda_\theta \lambda_e \bar{y}_j}{\lambda_\theta + m_j \lambda_e} + m_0 s_0\right]\right] + \frac{\lambda_e m_i \bar{y}_i}{\lambda_\theta + m_i \lambda_e}.$$

TABLE 6
*Constants used to construct total variation bounds*

| Sampler | $\gamma$ | $\phi$ | $\rho_1$ | $c_3$ | $d_R$ | $r$ | $\varepsilon$ |
|---|---|---|---|---|---|---|---|
| Block Gibbs | 0.3956 | 0.3589 | na | na | 28.328 | 0.0111 | 0.0246 |
| Gibbs | 0.41528 | na | 0.41527 | 2.6667 | 26.010 | 0.0009 | $5.6 \times 10^{-17}$ |

Observe that $E(\mu|\lambda)$ is a convex combination of $\bar{y}_i$ and $m_0$ and, furthermore, $E(\theta_i|\lambda)$ is a convex combination of $E(\mu|\lambda)$ and $\bar{y}_i$. If we let $\Delta$ denote the length of the convex hull of the set $\{\bar{y}_1, \bar{y}_2, \ldots, \bar{y}_K, m_0\}$, then for any $i = 1, 2, \ldots, K$, $[E(\theta_i|\lambda) - E(\mu|\lambda)]^2 \leq \Delta^2$ and $[E(\theta_i|\lambda) - \bar{y}_i]^2 \leq \Delta^2$.

## APPENDIX B

**B.1. Optimal starting values.** We desire the value of $(\theta, \mu)$ that minimizes

$$V_1(\theta, \mu) = \phi_1 v_1(\theta, \mu) + \phi_2 v_2(\theta) = \phi_1 \sum_{i=1}^{K} (\theta_i - \mu)^2 + \phi_2 \sum_{i=1}^{K} m_i(\theta_i - \bar{y}_i)^2.$$

Clearly, no matter what values are chosen for the $\theta_i$'s, the minimizing value of $\mu$ is $\bar{\theta}$. Thus, we need to find the value of $\theta$ that minimizes

$$\phi_1 \sum_{i=1}^{K} (\theta_i - \bar{\theta})^2 + \phi_2 \sum_{i=1}^{K} m_i(\theta_i - \bar{y}_i)^2.$$

Setting the derivative with respect to $\theta_i$ equal to 0 yields

$$(36) \qquad \theta_i = \frac{\phi_1 \bar{\theta} + \phi_2 m_i \bar{y}_i}{\phi_1 + \phi_2 m_i}.$$

Summing both sides over $i$ and dividing by $K$ yields an equation in $\bar{\theta}$ whose solution can be plugged back into (36) and this yields the optimal starting value

$$\theta_i = \frac{\phi_1 [\sum_{j=1}^{K} (m_j \bar{y}_j / (\phi_1 + \phi_2 m_j)) / \sum_{j=1}^{K} (m_j / (\phi_1 + \phi_2 m_j))] + \phi_2 m_i \bar{y}_i}{\phi_1 + \phi_2 m_i}.$$

## APPENDIX C

**C.1. Proof of Lemma 4.1.** Let

$$f_c(x) = \frac{(b + c/2)^\alpha}{\Gamma(\alpha)} x^{\alpha - 1} e^{-x(b + c/2)},$$

$$f_\beta(x) = \frac{(b + \beta/2)^\alpha}{\Gamma(\alpha)} x^{\alpha - 1} e^{-x(b + \beta/2)},$$

$$f_0(x) = \frac{b^\alpha}{\Gamma(\alpha)} x^{\alpha - 1} e^{-xb}.$$

Note that $x^*$ is the only positive solution to $f_c(x) = f_0(x)$. To prove the result it suffices to show that (i) $R_0(\beta) = f_\beta(x)/f_0(x) > 1$ for all $x \in (0, x^*)$ and all $\beta \in (0,c)$ and that (ii) $R_c(\beta) = f_\beta(x)/f_c(x) > 1$ for all $x \in (x^*, \infty)$ and all $\beta \in (0,c)$. Fix $k > 0$ and define a function

$$h(u) = \frac{ku}{1+ku} - \log(1+ku)$$

for $u \geq 0$. Since $h(0) = 0$ and $h'(u) < 0$, we know $h(u) < 0$ for $u \geq 0$. Hence,

(37) $$\frac{1}{u}\frac{k}{1+ku} - \frac{1}{u^2}\log(1+ku) < 0$$

for $u \geq 0$. Define another function,

$$g(u) = \frac{1}{u}\log(1+ku)$$

for $u > 0$. Since the the left-hand side of (37) is equal to $g'(u)$, we have established that $g(u)$ is decreasing for $u > 0$. Thus, if $x < x^* = \frac{2\alpha}{c}\log(1+\frac{c}{2b})$ and $\beta \in (0,c)$, then

$$\log R_0(\beta) = \alpha \log\left(1 + \frac{\beta}{2b}\right) - \frac{x\beta}{2}$$

$$> \alpha \log\left(1 + \frac{\beta}{2b}\right) - \frac{\alpha\beta}{c}\log\left(1 + \frac{c}{2b}\right)$$

$$= \alpha\beta\left[\frac{1}{\beta}\log\left(1 + \frac{\beta}{2b}\right) - \frac{1}{c}\log\left(1 + \frac{c}{2b}\right)\right] > 0,$$

and (i) is established. Case (ii) is similar.

## APPENDIX D

**D.1. Proof of Lemma 5.1.** First, let $g(v) = v + cv^{-1}$, where $c > 0$ and $v > 0$. It is easy to show that $g$ is minimized at $\hat{v} = \sqrt{c}$. Thus,

$$1 - \left(\frac{ax}{ax+y}\right)^2 - \left(\frac{y}{bx+y}\right)^2$$

$$= \frac{2bx^2y^2[y/x + ab(x/y)] + x^2y^2(b^2 + 4ab - a^2)}{(ax+y)^2(bx+y)^2}$$

$$\geq \frac{2bx^2y^2[2\sqrt{ab}] + x^2y^2(b^2 + 4ab - a^2)}{(ax+y)^2(bx+y)^2}$$

$$\geq \frac{x^2y^2(5b^2 + 4ab - a^2)}{(ax+y)^2(bx+y)^2}$$

$$= \frac{x^2y^2(5b-a)(b+a)}{(ax+y)^2(bx+y)^2}$$

$$> 0.$$

## APPENDIX E

**E.1.** $S_G \subset C_G = C_{G_1} \cap C_{G_2} \cap C_{G_3}$. First,

$$S_G = \{(\theta, \lambda) : V_3(\theta, \lambda) \leq d\}$$
$$= \left\{(\theta, \lambda) : e^{c_3 \lambda_\theta} + e^{c_3 \lambda_e} + \frac{\delta_7}{K \delta_1 \lambda_\theta} + v_3(\theta, \lambda) \leq d\right\}$$
$$\subset \left\{(\theta, \lambda) : e^{c_3 \lambda_\theta} \leq d, e^{c_3 \lambda_e} \leq d, \frac{\delta_7}{K \delta_1 \lambda_\theta} \leq d, v_3(\theta, \lambda) \leq d\right\}$$
$$= \left\{(\theta, \lambda) : \frac{\delta_7}{K \delta_1 d} \leq \lambda_\theta \leq \frac{\log d}{c_3}, 0 < \lambda_e \leq \frac{\log d}{c_3}, v_3(\theta, \lambda) \leq d\right\}.$$

As in the proof of Proposition 5.1, Jensen's inequality yields

$$\left(\frac{s_0 m_0 + K \lambda_\theta \bar{\theta}}{s_0 + K \lambda_\theta} - \bar{y}\right)^2 \leq \frac{s_0}{s_0 + K \lambda_\theta}(m_0 - \bar{y})^2 + \frac{K \lambda_\theta}{s_0 + K \lambda_\theta}(\bar{\theta} - \bar{y})^2$$
$$\leq (m_0 - \bar{y})^2 + v_3(\theta, \lambda),$$

and hence $S_G$ is contained in

$$C_G := \left\{(\theta, \lambda) : \frac{\delta_7}{K \delta_1 d} \leq \lambda_\theta \leq \frac{\log d}{c_3}, 0 < \lambda_e \leq \frac{\log d}{c_3}, \right.$$
$$\left.\left(\frac{s_0 m_0 + K \lambda_\theta \bar{\theta}}{s_0 + K \lambda_\theta} - \bar{y}\right)^2 \leq (m_0 - \bar{y})^2 + d\right\}.$$

Let $c_4 = \delta_7/(K \delta_1 d)$ and put $c_l$ and $c_u$ equal to $\bar{y} - \sqrt{(m_0 - \bar{y})^2 + d}$ and $\bar{y} + \sqrt{(m_0 - \bar{y})^2 + d}$, respectively. Note that $C_G = C_{G_1} \cap C_{G_2} \cap C_{G_3}$, where

$$C_{G_1} = \left\{(\theta, \lambda) : c_4 \leq \lambda_\theta \leq \frac{\log d}{c_3}\right\},$$
$$C_{G_2} = \left\{(\theta, \lambda) : 0 < \lambda_e \leq \frac{\log d}{c_3}\right\},$$
$$C_{G_3} = \left\{(\theta, \lambda) : c_l \leq \frac{s_0 m_0 + K \lambda_\theta \bar{\theta}}{s_0 + K \lambda_\theta} \leq c_u\right\}.$$

## APPENDIX F

**F.1. Closed form expression for $\varepsilon_{\mathbf{G}}$.** Recall that

$$\varepsilon_G = \left[\frac{s_0 + K c_4}{s_0 + K \log(d)/c_3}\right]^{1/2} \left[\int_{\mathbb{R}} \int_{\mathbb{R}^K} g_1(\mu, \theta) g_2(\mu) \, d\theta \, d\mu\right].$$

A straightforward calculation shows that

$$\int_{\mathbb{R}^K} g_1(\mu,\theta)\,d\theta$$
(38)
$$= \left(\frac{c_4 c_3}{\log d}\right)^{K/2} \prod_{i=1}^{K} \sqrt{\frac{1}{1+m_i}} \exp\left\{-\frac{m_i \log d}{2c_3(1+m_i)}(\mu - \bar{y}_i)^2\right\}.$$

Thus,

$$\int_{\mathbb{R}^K}\int_{\mathbb{R}} g_1(\mu,\theta) g_2(\mu)\,d\theta\,d\mu$$
$$= \left(\frac{c_4 c_3}{\log d}\right)^{K/2} \prod_{i=1}^{K} \sqrt{\frac{1}{1+m_i}} \int_{\mathbb{R}} g_2(\mu) \prod_{i=1}^{K} \exp\left\{-\frac{m_i \log d}{2c_3(1+m_i)}(\mu - \bar{y}_i)^2\right\} d\mu.$$

Now

$$\int_{\mathbb{R}} g_2(\mu) \exp\left\{-\frac{\log d}{2c_3}\sum_{i=1}^{K}\frac{m_i}{1+m_i}(\mu-\bar{y}_i)^2\right\} d\mu$$

$$= \sqrt{\frac{s_0 + K\log(d)/c_3}{2\pi}}$$

$$\times \Bigg[\int_{-\infty}^{\bar{y}} \exp\left\{-\frac{s_0+K\log(d)/c_3}{2}(\mu-c_u)^2\right\}$$

$$\times \exp\left\{-\frac{\log d}{2c_3}\sum_{i=1}^{K}\frac{m_i}{1+m_i}(\mu-\bar{y}_i)^2\right\} d\mu$$

$$+ \int_{\bar{y}}^{\infty} \exp\left\{-\frac{s_0+K\log(d)/c_3}{2}(\mu-c_l)^2\right\}$$

$$\times \exp\left\{-\frac{\log d}{2c_3}\sum_{i=1}^{K}\frac{m_i}{1+m_i}(\mu-\bar{y}_i)^2\right\} d\mu \Bigg].$$

Define

$$v = \left[s_0 + \frac{\log d}{c_3}\left(K + \sum_{i=1}^{K}\frac{m_i}{1+m_i}\right)\right]^{-1}$$

and put

$$m_l = v\left[c_l s_0 + \frac{\log d}{c_3}\left(Kc_l + \sum_{i=1}^{K}\frac{\bar{y}_i m_i}{1+m_i}\right)\right]$$

and

$$m_u = v\left[c_u s_0 + \frac{\log d}{c_3}\left(Kc_u + \sum_{i=1}^{K}\frac{\bar{y}_i m_i}{1+m_i}\right)\right].$$

Then

$$\int_{-\infty}^{\bar{y}} \exp\left\{-\frac{s_0 + K\log(d)/c_3}{2}(\mu - c_u)^2\right\} \exp\left\{-\frac{\log d}{2c_3}\sum_{i=1}^{K}\frac{m_i}{1+m_i}(\mu - \bar{y}_i)^2\right\} d\mu$$

$$= \exp\left\{-\frac{c_u^2 s_0}{2} - \frac{\log d}{2c_3}\left[Kc_u^2 + \sum_{i=1}^{K}\frac{\bar{y}_i^2 m_i}{1+m_i}\right] + \frac{m_u^2}{2v}\right\}\sqrt{2\pi v}\,\Phi\left(\frac{\bar{y} - m_u}{\sqrt{v}}\right)$$

and

$$\int_{\bar{y}}^{\infty} \exp\left\{-\frac{s_0 + K\log(d)/c_3}{2}(\mu - c_l)^2\right\} \exp\left\{-\frac{\log d}{2c_3}\sum_{i=1}^{K}\frac{m_i}{1+m_i}(\mu - \bar{y}_i)^2\right\} d\mu$$

$$= \exp\left\{-\frac{c_l^2 s_0}{2} - \frac{\log d}{2c_3}\left[Kc_l^2 + \sum_{i=1}^{K}\frac{\bar{y}_i^2 m_i}{1+m_i}\right] + \frac{m_l^2}{2v}\right\}$$

$$\times \sqrt{2\pi v}\left(1 - \Phi\left(\frac{\bar{y} - m_l}{\sqrt{v}}\right)\right).$$

Putting all of this together yields

$$\varepsilon_G = \sqrt{v(s_0 + Kc_4)}\sqrt{\prod_{i=1}^{K}\frac{1}{1+m_i}}\left(\frac{c_4 c_3}{\log d}\right)^{K/2}\exp\left\{-\frac{\log d}{2c_3}\sum_{i=1}^{K}\frac{\bar{y}_i^2 m_i}{1+m_i}\right\}$$

$$\times \left[\exp\left\{-\frac{c_u^2 s_0}{2} - \frac{Kc_u^2 \log d}{2c_3} + \frac{m_u^2}{2v}\right\}\Phi\left(\frac{\bar{y} - m_u}{\sqrt{v}}\right)\right.$$

$$\left.+ \exp\left\{-\frac{c_l^2 s_0}{2} - \frac{Kc_l^2 \log d}{2c_3} + \frac{m_l^2}{2v}\right\}\left(1 - \Phi\left(\frac{\bar{y} - m_l}{\sqrt{v}}\right)\right)\right],$$

where $\Phi(\cdot)$ denotes the standard normal cumulative distribution function.

**Acknowledgment.** The authors are grateful to two anonymous referees whose insightful comments led to substantial improvements in the paper.

School of Statistics  
University of Minnesota  
Twin Cities Campus  
313 Ford Hall  
224 Church Street, SE  
Minneapolis, Minnesota 55455  
USA  
e-mail: galin@stat.umn.edu

Department of Statistics  
University of Florida  
Gainesville, Florida 32611  
USA  
e-mail: jhobert@stat.ufl.edu